\documentstyle[12pt,leqno,amssymb,graphics]{amsart}

\newtheorem{thm}{Theorem}[section]
\newtheorem{lemma}[thm]{Lemma}
\newtheorem{prop}[thm]{Proposition}
\newtheorem{cor}[thm]{Corollary}

\theoremstyle{definition}
\newtheorem{dfn}[thm]{Definition}
\theoremstyle{remark}

\begin{document}

\newcommand{\pr}{\protect\ref}
\newcommand{\su}{\subseteq}
\newcommand{\pa}{{\partial}}

\newcommand{\R}{{\Bbb R}}
\newcommand{\E}{{{\Bbb R}^3}}
\newcommand{\C}{{{\Bbb Z}/2}}

\newcommand{\r}{{\mathrm{rank}}}
\newcommand{\p}{{\psi}}
\newcommand{\pp}{{\widehat{\psi}}}
\newcommand{\hp}{\widehat{\Psi}}
\newcommand{\br}{\bar}

\newcommand{\hh}{{H_1(F,\C)}}
\newcommand{\ov}{{O(V,g)}}
\newcommand{\oh}{{O(\hh,g)}}
\newcommand{\ohi}{{O(\hh,g^i)}}

\newcommand{\ep}{{\epsilon}}
\newcommand{\eep}{{\widehat{\epsilon}}}

\newcounter{numb}

\title[Embeddings of Surfaces and Quadruple Points]
{Embeddings of Surfaces into 3-Space \\ and
Quadruple Points of Regular Homotopies}
\author{Tahl Nowik}
\address{Department of Mathematics, Columbia University, New York, 
NY 10027, USA.} 
\email{tahl@@math.columbia.edu}
\date{November 9, 1999}

\begin{abstract}
Let $F$ be a closed orientable surface. 
We give an explicit formula for the number
mod 2 of quadruple points occurring in any generic regular homotopy between
any two regularly homotopic embeddings $e,e':F\to\E$.
The formula is in terms of homological data extracted 
from the two embeddings.
\end{abstract}

\maketitle

\section{Introduction}\label{A}
For $F$ a closed surface and $i,i':F\to\E$ two 
regularly homotopic generic immersions, we are interested in the number mod 2
of quadruple points occurring in generic regular homotopies 
between $i$ and $i'$. 
It has been shown in [N1] that this number is the same for all such
regular homotopies, and so it is a function of $i$ and 
$i'$ which we denote $Q(i,i')\in\C$. 
For $F$ orientable and $e,e':F\to\E$ two regularly homotopic \emph{embeddings},
we give an explicit formula for $Q(e,e')$ 
which depends on the following 
data: If $e:F\to\E$ is an embedding then $e(F)$ splits $\E$ into 
two pieces, one compact which will be denoted $M^0(e)$ and the other 
non-compact which will be denoted $M^1(e)$. By restriction of range 
$e$ induces maps
$e^k:F\to M^k(e)$ ($k=0,1$) and let $A^k(e) \su \hh$ be the kernel of the map
induced by $e^k$ on $H_1(\cdot,\C)$. Let $o(e)$ be the orientation
on $F$ which is induced from $M^0(e)$ to $\pa M^0(e) = e(F)$ 
and then via $e$ to $F$. 
Our formula for $Q(e,e')$ will be in terms of the two triplets 
$A^0(e),A^1(e),o(e)$ and $A^0(e'),A^1(e'),o(e')$.
Our formula will be also easily extended to finite unions of 
closed orientable surfaces.

For two special cases a formula for $Q(e,e')$ (for $e,e'$ embeddings)
has already been known:
The case where $F$ is a sphere has appeared in [MB] and [N1], and the 
case where $F$ is a torus has appeared in [N1].
The starting point for our work will be [N2] where an
explicit formula has been given for $Q(i,i\circ h)$, where $i:F\to\E$ is any
generic immersion and $h:F\to F$ is any diffeomorphism such that 
$i$ and $i\circ h$ are regularly homotopic.

\section{Totally Singular Decompositions}\label{B}

Let $V$ be a finite dimensional vector space over $\C$. 
A function $g:V\to\C$ is called a \emph{quadratic form} if $g$ satisfies:
$g(x+y)=g(x)+g(y)+B(x,y)$ for all $x,y \in V$, 
where $B(x,y)$ is a bilinear form.
The following properties follow:
(a) $g(0)=0$. (b) $B(x,x)=0$ for all $x\in V$.
(c) $B(x,y)=B(y,x)$ for all $x,y\in V$.  
$g$ is called \emph{non-degenerate} if $B$ is non-degenerate,
i.e. for any $0\neq x\in V$ there is $y\in V$ with $B(x,y)\neq 0$.
For an exposition of quadratic forms see [C].

In what follows we always assume that our vector space $V$ is equipped with a
non-degenerate quadratic form $g$. It then follows that $\dim V$ is even.
A subspace $A\su V$ such that $g|_A \equiv 0$ is called a
\emph{totally singular} subspace.
A pair $(A,B)$ of subspaces of $V$ will be called a 
\emph{totally singular decomposition} (abbreviated TSD) of $V$ if 
$V=A\oplus B$ and both
$A$ and $B$ are totally singular. 
It then follows that $\dim A = \dim B$. 
(We remark that TSDs do not always exist. They will however always exist
for the quadratic forms which will arise in our geometric considerations,
as seen in Lemma \pr{g1} below.)
A linear map $T:V\to V$ is called \emph{orthogonal}
if $g(T(x))=g(x)$ for all $x\in V$.
It then follows that $B(T(x),T(y))=B(x,y)$ for all $x,y\in V$
and that $T$ is invertible.
The group of all orthogonal maps of $V$ with respect to $g$ will
be denoted $\ov$.

The proof of the following lemma appears in [C]:

\begin{lemma}\label{ll1}
Let $\dim V =2n$.
\begin{enumerate}
\item If $A\su V$ is a totally singular subspace of dimension $n$
then there exists a $B\su V$ such that $(A,B)$ is a TSD of $V$.
\item  If $(A,B)$ is a TSD of $V$ and $a_1,\dots,a_n$ is a given basis
for $A$ then there is a basis $b_1,\dots,b_n$ for $B$ such that
$B(a_i,b_j)=\delta_{ij}$. 
\end{enumerate}
\end{lemma}

\begin{dfn}\label{d1}
If $(A,B)$ is a TSD of $V$ then a basis $a_1,\dots,a_n,b_1,\dots,b_n$ 
of $V$ will be called 
\emph{$(A,B)$-good} if $a_i\in A$, $b_i\in B$ and
$B(a_i,b_j)=\delta_{ij}$.
\end{dfn}

The following two lemmas follow 
directly from the definition of quadratic form: 

\begin{lemma}\label{ll2}
Let $(A,B)$ be a TSD of $V$ and
$a_1,\dots,a_n,b_1,\dots,b_n$ an $(A,B)$-good basis for $V$.
If $v=\sum x_i a_i + \sum y_i b_i$ and
$v'=\sum x'_i a_i + \sum y'_i b_i$ 
then $g(v) = \sum x_i y_i$ and $B(v,v')=\sum x_i y'_i + \sum y_i x'_i$.
\end{lemma}

\begin{lemma}\label{ll3}
Let $(A,B)$ and $(A',B')$ be two TSDs of $V$. 
Let $a_1,\dots,a_n,b_1,\dots,b_n$ be an $(A,B)$-good basis for $V$
and $a'_1,\dots,a'_n,b'_1,\dots,b'_n$ an $(A',B')$-good basis for $V$.
If $T:V\to V$ is the linear map
defined by $a_i\mapsto a'_i$, $b_i\mapsto b'_i$ then $T\in\ov$.
\end{lemma}

For $T\in\ov$ we define $\p(T)\in\C$ by: 
$$\p(T) \ = \ \r(T-Id)\mod{2}.$$
It has been shown in [N2] that $\p:\ov\to\C$ is a (non-trivial) 
homomorphism.

\begin{lemma}\label{l1}
If $(A,B)$ is a TSD of $V$ and
$T\in\ov$ satisfies $T(A)=A$ and $T(B)=B$ then $\p(T)=0$.
\end{lemma}

\begin{pf}
By Lemma \pr{ll1} there exists an $(A,B)$-good basis 
$a_1,\dots,a_n,b_1,\dots,b_n$ for $V$. 
Using Lemma
\pr{ll2} it is easy to verify 
that the matrix of $T$
with respect to such a basis has the form:
$\left(
\begin{matrix}
S^t & 0 \\
0 & S^{-1}
\end{matrix}
\right)$ where $S\in GL_n(\C)$.
It follows that $\p(T)=0$.
\end{pf}

Given two TSDs $(A,B)$, $(A',B')$ of $V$ then by 
Lemmas \pr{ll1} and \pr{ll3} there exists a $T\in\ov$ such
that $T(A)=A'$ and $T(B)=B'$. It follows from Lemma
\pr{l1} that if $T_1,T_2$ are two such $T$s then $\p(T_1)=\p(T_2)$.
And so the following is well defined:

\begin{dfn}\label{d2}
For a  pair $(A,B)$, $(A',B')$ of TSDs of $V$
let $\pp(A,B;A',B')\in\C$ be defined by
$\pp(A,B;A',B')=\p(T)$ for some (thus all) $T\in\ov$ 
with $T(A)=A'$ and $T(B)=B'$.
\end{dfn}

\begin{dfn}\label{d3}
For two TSDs $(A,B)$, $(A',B')$ of $V$, we will write
$(A,B) \sim (A',B')$ if $\pp(A,B; A',B')=0$. 
\end{dfn}

Since $\p$ is a homomorphism, 
$\pp(A,B; A'',B'')=\pp(A,B; A',B')+\pp(A',B'; A'',B'')$
for any three TSDs $(A,B), (A',B'), (A'',B'')$.
It follows that $\sim$ is an equivalence relation 
with precisely two equivalence classes and that
$\pp(A,B; A'',B'')=\pp(A',B'; A'',B'')$ whenever $(A,B) \sim (A',B')$.

\begin{lemma}\label{l2}
Let $\dim V = 2n$ and let $A\su V$ be a totally singular subspace of 
dimension $n$. If $T\in\ov$ satisfies
$T(x)=x$ for every $x\in A$ then $\p(T)=0$.
\end{lemma}

\begin{pf}
By Lemma \pr{ll1} there is a $B\su V$ such that $(A,B)$ is a 
TSD of $V$ and an $(A,B)$-good basis $a_1,\dots,a_n,b_1,\dots,b_n$
for $V$. Using Lemma \pr{ll2} it is easy to verify that
the matrix of $T$
with respect to such a basis has the form:
$\left(
\begin{matrix}
I & S \\
0 & I
\end{matrix}\right)$
where $I$ is the $n\times n$ identity matrix and $S\in M_n(\C)$
is an alternating matrix, i.e. if $S=\{s_{ij}\}$ then
$s_{ii}=0$ and $s_{ij} = s_{ji}$. 
Since alternating matrices have even rank, it follows that $\p(T)=0$.
\end{pf}

\begin{cor}\label{p1}
Let $(A,B)$ and $(A',B')$ be two TSDs of $V$.
If $A=A'$ or $B=B'$ then $(A,B)\sim (A',B')$.
\end{cor}

\begin{pf}
Say $A=A'$. By Lemmas \pr{ll1} and \pr{ll3}
there exists a $T\in\ov$ with $T(x)=x$ for all $x\in A=A'$ and 
$T(B)=B'$. The conclusion follows from Lemma \pr{l2}.
\end{pf}

Let $V_0,V_1\su V$ be two subspaces of $V$. We will write $V_0 \bot V_1$
if $B(x,y)=0$ for every $x\in V_0$, $y\in V_1$. 
The following is clear:

\begin{lemma}\label{p2}

Let $V_0,V_1 \su V$ satisfy 
$V=V_0\oplus V_1$ and $V_0 \bot V_1$.
\begin{enumerate}
\item If for $l=0,1$, $(A_l,B_l)$ is a TSD of $V_l$ (with respect to
$g|_{V_l}$ which is indeed non-degenerate) then
$(A_0 + A_1, B_0 + B_1)$ is a TSD of $V$.
\item If $(A'_l,B'_l)$ is another TSD of $V_l$
and $(A_l,B_l) \sim (A'_l,B'_l)$ ($l=0,1$) then
$(A_0 + A_1, B_0 + B_1) \sim (A'_0 + A'_1, B'_0 + B'_1)$.
\end{enumerate}
\end{lemma}

\section{Statement of Main Result}\label{C}

A \emph{surface} is by definition assumed connected. A finite union of
surfaces will be called a \emph{system of surfaces}. 
Let $S$ be a system of closed 
surfaces and $H_t:S\to\E$ a generic regular homotopy. 
We denote by $q(H_t)\in\C$
the number mod 2 of quadruple points occurring in $H_t$.
The following has been shown in [N1]:

\begin{thm}\label{tn1}
Let $S$ be a system of closed surfaces (not necessarily orientable.)
If $H_t,G_t:S\to\E$ are two generic regular homotopies 
between the same two generic immersions, then $q(H_t)=q(G_t)$.
\end{thm}

\begin{dfn}\label{qii}
Let $S$ be a system of closed surfaces and
$i,i':S\to\E$ two regularly homotopic generic immersions.
We define $Q(i,i')\in\C$ by $Q(i,i')=q(H_t)$, where $H_t$ is any 
generic regular homotopy between $i$ and $i'$. 
This is well defined by Theorem \pr{tn1}. 
\end{dfn}

Let $F$ from now on denote a closed orientable surface.
A simple closed curve in $F$ will be called a 
\emph{circle}.
If $c$ is a circle
in $F$, the homology class of $c$ in $\hh$ will be denoted by $[c]$.
Any immersion $i:F\to \E$ induces a quadratic form $g^i:\hh\to\C$
whose associated bilinear form $B(x,y)$ is the algebraic intersection form 
$x \cdot y$ of $\hh$, as follows: For $x\in\hh$ let $A\su F$ be an 
annulus bounded by circles $c,c'$ with $[c]=x$,
let $j:A\to\E$ be an embedding which is regularly homotopic to $i|_A$ and
define $g^i(x)$ to be the $\C$ linking number between $j(c)$ and $j(c')$
in $\E$.
One needs to verify that $g^i(x)$ is independent of the choices 
being made and that $g^i(x+y)=g^i(x)+g^i(y)+x\cdot y$.
This has been done in [P].
Also, $i,i':F\to\E$ are regularly homotopic iff $g^i=g^{i'}$.

If $e:F\to\E$ is an embedding then $e(F)$ splits $\E$ into two pieces
one compact and one non-compact. We denote the compact piece by
$M^0(e)$ and the non-compact piece by $M^1(e)$.
By restriction of range, $e$ induces maps 
$e^k : F \to M^k(e)$, $k=0,1$. 
Let $e^k_* : \hh \to H_1(M^k(e),\C)$ be the maps 
induced on homology
and finally let $A^k(e) = \ker e^k_*$, $k=0,1$.

\begin{lemma}\label{g1}
Let $e:F\to \E$ be an embedding, then $(A^0(e), A^1(e))$ is a TSD of $\hh$
with respect to the quadratic form $g^e$.

\end{lemma}

\begin{pf}
We first show that each $A^k(e)$ is totally singular: For
$x\in A^k(e)$ let $A,c,c'$ be as in the definition of $g^e(x)$
and simply take $j=e|_A$.
Since $e^k_*(x)=0$, $e(c)$ bounds a properly embedded (perhaps non-orientable)
surface $S$ in $M^k(e)$. Since $e(c')$ is disjoint from $S$, the $\C$ linking
number between $e(c)$ and $e(c')$ in $\E$ is 0, and so $g^e(x)=0$.
Now, the fact that $\hh=A^0(e) \oplus A^1(e)$ is a consequence of
the $\C$ Mayer-Vietoris sequence for $\E = M^0(e) \cup M^1(e)$
where $F$ is identified with $M^0(e) \cap M^1(e)$ via $e$.
\end{pf}

If $e,e':F\to\E$ are two regularly homotopic embeddings
then $g^e=g^{e'}$ so $(A^0(e),A^1(e))$ and  $(A^0(e'),A^1(e'))$ 
are TSDs of $\hh$ with respect to the same quadratic form and so
$\pp(  A^0(e),A^1(e)  ; A^0(e'),A^1(e') )$ is defined.
We spell out the actual computation involved in 
$\pp(  A^0(e),A^1(e)  ; A^0(e'),A^1(e') )$: 
\begin{enumerate}
\item Find a basis 
$a_1,\dots,a_n,b_1,\dots,b_n$ for $\hh$ such that 
$e^0_*(a_i)=0$, $e^1_*(b_i)=0$ 
and $a_i \cdot b_j = \delta_{ij}$. 
\item Find a similar basis
$a'_1,\dots,a'_n,b'_1,\dots,b'_n$ using $e'$ in place of $e$.
\item Let $m$ be the dimension of the subspace of $\hh$ spanned by:
$$a'_1 - a_1 \ , \ \dots \ ,  \ a'_n - a_n  \ ,  \ b'_1 - b_1  \ ,   \ 
\dots \ ,  \ b'_n - b_n.$$
\item $\pp(  A^0(e),A^1(e)  ; A^0(e'),A^1(e') ) = m\bmod{2}$, 
(an element in $\C$.)
\end{enumerate}

\begin{dfn}\label{d6}
If $e:F\to \E$ is an embedding then we define $o(e)$ to be the 
orientation on $F$ which is induced from $M^0(e)$ to $\pa M^0(e) = e(F)$ 
and then via $e$ to $F$
(and where the orientation on $M^0(e)$ is the restriction of the orientation
of $\E$.)
If $e,e':F\to \E$ are two embeddings then we define $\eep(e,e')\in\C$
to be 0 if $o(e) = o(e')$ and 1 if $o(e) \neq o(e')$.
\end{dfn}

Our purpose in this work is to show:

\begin{thm}\label{t1}
Let $n$ be the genus of $F$. If $e,e':F\to \E$ are 
two regularly homotopic embeddings then:
$$Q(e,e') \ = \ \pp(  A^0(e),A^1(e)  ;  A^0(e'),A^1(e')  ) \ 
+ \ (n+1)\eep(e,e').$$
\end{thm}

Our starting point is the following theorem which has been proved
in [N2]:

\begin{thm}\label{tn2}
For any generic immersion $i:F\to\E$ 
and any diffeomorphism
$h:F\to F$ such that $i$ and $i\circ h$ are regularly homotopic,
$$Q(i,i\circ h) = \p(h_*) + (n+1)\ep(h),$$
where $h_*$ is the map induced by $h$ on $H_1(F,\C)$, $n$ is the genus
of $F$ and $\ep(h)\in\C$ is 0 or 1 according to whether $h$ is orientation
preserving or reversing, respectively.
\end{thm}

\section{Equivalent Embeddings and $k$-Extendible Regular Homotopies}\label{D}

Let $e:F\to\E$ be an embedding, let $P\su\E$ be a plane and
assume $e(F)$ intersects $P$ transversally in a unique circle.
Let $c=e^{-1}(P)$ then $c$ is a separating circle in $F$.
Let $A$ be a regular neighborhood of $c$ in $F$ and let $F_0,F_1$ be
the connected components of $F-int A$. 
(A lower index will always be related to the splitting
of $\E$ via a plane, the assignment of 0 and 1 to the two sides
being arbitrary. An upper index on the other hand is related to the
splitting of $\E$ via the image of a closed surface,
assigning 0 to the compact side and 1 to the non-compact side.)
Let $\br{F_l}$ ($l=0,1$) be the
closed surface obtained by gluing a disc $D_l$ to $F_l$. Let 
$e_l:\br{F_l}\to\E$
be the embedding such that $e_l|_{F_l}=e|_{F_l}$ and $e_l(D_l)$ is parallel
to $P$. Let $i_{F_l F}:F_l \to F$ and $i_{F_l \br{F_l}}: F_l \to\br{F_l}$
denote the inclusion maps.
The induced map ${i_{F_l\br{F_l}}}_*:H_1(F_l,\C)\to H_1(\br{F_l},\C)$ 
is an isomorphism and let $h_l:H_1(\br{F_l},\C)\to \hh$ be the map
$h_l={i_{F_l F}}_* \circ ({i_{F_l \br{F_l}}}_*)^{-1}$.

\begin{lemma}\label{sp}
Under the above assumptions and definitions:
$A^k(e)=h_0(A^k(e_0)) + h_1(A^k(e_1))$, $k=0,1$.
\end{lemma}

\begin{pf}
This follows from the fact that the inclusions
$F_0 \cup F_1 \to \br{F_0} \cup \br{F_1}$, \ 
$F_0 \cup F_1 \to F$, \ 
$M^0(e_0) \cup M^0(e_1) \to M^0(e)$ and
$M^1(e) \to \E-(M^0(e_0) \cup M^0(e_1))$
all induce isomorphisms on $H_1(\cdot,\C)$
and the splitting of each of the above spaces via $P$ 
induces a direct sum decomposition. 
We only check that the inclusion $M^1(e) \to \E-(M^0(e_0) \cup M^0(e_1))$
induces isomorphism on $H_1(\cdot,\C)$. 
Indeed $\E-(M^0(e_0) \cup M^0(e_1))$ is obtained from
$M^1(e)$ by gluing a 2-handle along $e(A)$, and the inclusion of $e(A)$
in $M^1(e)$ is null-homotopic.
\end{pf}

\begin{dfn}\label{eq}
Two embeddings $e,f:F\to\E$ will be called \emph{equivalent} if:
\begin{enumerate}
\item There is a regular homotopy between $e$ and $f$ with no
quadruple points. 
\item $(A^0(e),A^1(e)) \sim (A^0(f),A^1(f))$.
\item $o(e) = o(f)$
\end{enumerate}
\end{dfn}

\begin{dfn}\label{st}
An embedding $e:F\to\E$ will be called \emph{standard} if its image $e(F)$ is 
a surface in $\E$ as in Figure \pr{f1}. 
\end{dfn}

\begin{figure}[h]
\scalebox{0.6}{\includegraphics{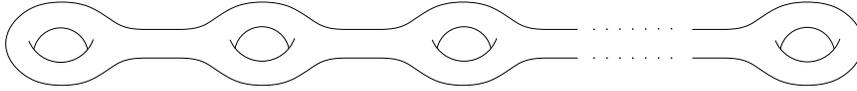}}
\caption{Image of a standard embedding.}\label{f1}
\end{figure}

In Proposition \pr{p3} below we will show that any embedding $e:F\to\E$
is equivalent to a standard embedding. The following lemma will be used in the
induction step:

\begin{lemma}\label{s2}
Let $e:F\to\E$ be an embedding.
Assume $e(F)$ intersects a plane $P\su\E$ transversally in one circle
and let $c,A,F_l,\br{F_l},D_l,e_l$ be as above.
If $e_l:\br{F_l}\to \E$ ($l=0,1$) are both equivalent to standard 
embeddings, then $e$ is equivalent to a standard embedding.
\end{lemma}

\begin{pf}
Changing $e$ by isotopy, we may assume $e(A)$ is a very thin tube. 
$e_l:\br{F_l}\to \E$ is equivalent to a standard embedding $f_l$ via
a regular homotopy $(H_l)_t:\br{F_l}\to \E$ having no quadruple points.
We may further assume that each $(H_l)_t$ moves $\br{F_l}$
only within the corresponding half-space defined by $P$, that
each $f_l(D_l)$ is situated
at the point of $f_l(\br{F_l})$ which is closest to $P$ and that these
two points are opposite each other with respect to $P$.
We now perform both $(H_l)_t$, letting the thin
tube $A$ be carried along. If we make sure the thin tube $A$ does not
pass triple points occurring in $F_1$ and $F_2$ then the regular homotopy
$H_t$ induced on $F$ in this way will also have no quadruple points.
Since $e(A)$ has approached $e_l(\br{F_l})$ from $M^1(e_l)$ and
since $o_{e_l} = o_{f_l}$, we also have at the end of $H_t$ 
that $A$ approaches $f_l(\br{F_l})$ from $M^1(f_l)$.
And so we may continue moving
the tube $A$ until it is all situated in the region between 
$f_0(\br{F_0})$ and $f_1(\br{F_1})$,
then canceling all knotting by
having the thin tube pass itself (this involves only double lines)
until $A$ is embedded as a straight tube connecting $f_0(\br{F_0})$
to $f_1(\br{F_1})$ and so
the final map $f:F\to\E$ thus obtained is indeed a standard embedding. 
By assumption $(A^0(e_l), A^1(e_l)) \sim (A^0(f_l) , A^1(f_l))$, $l=0,1$
which implies that
$(h_l(A^0(e_l)), h_l(A^1(e_l))) \sim (h_l(A^0(f_l)) , h_l(A^1(f_l)))$, $l=0,1$
as TSDs of $V_l = h_l(H_1(\br{F_l},\C))\su \hh$.
(Note that $h_l$ preserves the corresponding quadratic forms.)
But $\hh=V_0\oplus V_1$ and $V_0 \bot V_1$ and so by
Lemma \pr{p2} and Lemma \pr{sp} 
$(A^0(e),A^1(e)) \sim (A^0(f),A^1(f))$.
Finally, from $o_{e_l}=o_{f_l}$ it follows that $o(e) = o(f)$. 
\end{pf}

\begin{dfn}\label{d5}
Let $e,f:F\to\E$ be two embeddings.
A regular homotopy $H_t:F\to\E$ ($a\leq t \leq b$) 
with $H_a=e$, $H_b=f$ will be called
\emph{k-extendible} (where $k$ is either 0 or 1) if there exists a 
regular homotopy $G_t: M^k(e)\to\E$ ($a\leq t \leq b$) satisfying: 

\begin{enumerate}
\item $G_a$ is the inclusion map of $M^k(e)$ in $\E$.
\item $H_t = G_t \circ e^k$. 
(Recall that $e^k:F\to M^k(e)$ is simply $e$ with range 
restricted to $M^k(e)$.)
\item $G_b$ is an embedding with $G_b(M^k(e)) = M^k(f)$. 
\end{enumerate}
\end{dfn}

\begin{lemma}\label{g2}
If for a given $k$
there is a $k$-extendible regular homotopy between the embeddings
$e$ and $f$ then $A^k(e)=A^k(f)$.
\end{lemma}

\begin{pf}
$f = H_b = G_b \circ e^k$ and so $f^k = G_b^k \circ e^k$
where $G_b^k : M^k(e) \to M^k(f)$ is the map $G_b$ with range restricted
to $M^k(f)$. Since $G_b^k$ is a diffeomorphism it
follows that $\ker f^k_* = \ker e^k_*$.
\end{pf}

\begin{cor}\label{g3}
If there is a $k$-extendible regular homotopy between the embeddings
$e$ and $f$ for either $k=0$ or $k=1$ then:
\begin{enumerate}
\item\label{1} $(A^0(e),A^1(e)) \sim (A^0(f),A^1(f))$.
\item\label{2} $o(e)=o(f)$.
\end{enumerate}
\end{cor}

\begin{pf}
\pr{1} follows from Lemma \pr{g2} and Corollary \pr{p1}. 
Since $G_a$ is the inclusion and $G_t$ is a regular homotopy it follows 
that $G_b$ is orientation preserving. This implies \pr{2}.
\end{pf}

\begin{prop}\label{p3}
Every embedding $e:F\to \E$ is equivalent to a standard embedding.
\end{prop}

\begin{pf}
The proof is by induction on the genus of $F$. If $F=S^2$ then any $e$ 
is isotopic to a standard embedding and isotopic embeddings
are equivalent. 
So assume $F$ is of positive genus and so there is a compressing disc
$D$ for $e(F)$ in $\E$ (i.e. $D\cap e(F)=\pa D$ and $\pa D$ does not
bound a disc in $e(F)$.)
Let $c=e^{-1}(\pa D) \su F$ and let 
$A$ be a regular neighborhood of $c$ in $F$. Isotoping $A$ along $D$ 
as before we may assume $A$ is embedded as a thin tube.
There are four cases to be considered according to whether
$D$ is contained in $M^0(e)$ or $M^1(e)$ and whether $\pa D$ separates
or does not separate $e(F)$.

\emph{Case 1}: $D\su M^0(e)$ and $\pa D$ separates
$e(F)$. It then follows that $D$ separates $M^0(e)$.
If $F_0,F_1$ denote the two components of $F-int A$ and 
$e_l:\br{F_l}\to\E$ are defined as before then it follows from the
assumptions of this case that $M^0(e_0)$ and $M^0(e_1)$
are disjoint and the tube $e(A)$ approaches each
$e_l(\br{F_l})$ from its non-compact side, i.e. from $M^1(e_l)$.
Move each foot of the tube $e(A)$ (see Figure \pr{f2})
along the corresponding surface $e_l(\br{F_l})$ until they are each situated 
at the point $p_l$ of $e_l(\br{F_l})$ having maximal $z$-coordinate. 
In particular it follows that now $e(A)$ approaches each $e_l(\br{F_l})$ 
from above. We now uniformly
shrink each $e(F_l)$ towards the point 
$p_l$ until it is contained in a 
tiny ball $B_l$ attached from below to the corresponding foot of $e(A)$,
arriving at a new embedding $e':F\to\E$.
This regular homotopy is clearly 0-extendible, and since no self 
intersection may occur within each of $F_0$, $F_1$ and $A$, this 
regular homotopy has no quadruple points. 
And so by Corollary \pr{g3} $e'$ is equivalent to $e$.
We now continue by isotopy, deforming the thin tube $e'(A)$ until it is
a straight tube, and rigidly carrying $B_0$ and $B_1$ along.
We finally arrive at an embedding $e''$ for which there is a plane $P$
intersecting $e''(F)$ as in Lemma \pr{s2} with our $F_0$ and $F_1$ on the 
two sides of $P$. Since the genus of both $\br{F_0}$
and $\br{F_1}$ is smaller than that of $F$ then by
our induction hypothesis and Lemma \pr{s2}, $e''$ is equivalent 
to a standard embedding.

\begin{figure}[h]
\scalebox{0.6}{\includegraphics{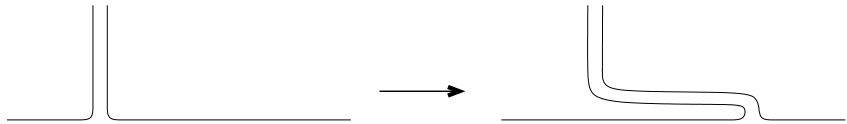}}
\caption{Moving the foot of a tube.}\label{f2}
\end{figure}

\emph{Case 2}: $D\su M^1(e)$ and $\pa D$ separates $e(F)$. This time 
either $M^0(e_0) \su M^0(e_1)$ or $M^0(e_1)\su M^0(e_0)$ and assume the former
holds. In this case $e(A)$ approaches only $e_0(\br{F_0})$ 
from its non-compact side and so we push the tube and perform 
the uniform shrinking as above only with $F_0$. This is a 1-extendible
regular homotopy since we are shrinking $M^0(e_0)$ which is part of $M^1(e)$.
Now, if $B$ is the tiny
ball into which we have shrunken $e(F_0)$ then $\pa B$ supplies separating
compressing discs on both sides of $e(F)$ and so we are done by Case 1.

\emph{Case 3}: $D\su M^0(e)$ and $\pa D$ does not separate $e(F)$.
If $F'=F-int A$ and $e':\br{F'}\to \E$ is induced as above
(where $\br{F'}$ is the surface obtained from $F'$ by gluing two discs to it)
then both feet of the tube $e(A)$ approach $e'(\br{F'})$ from its
non-compact side. Push the feet of $e(A)$ until they are both 
situated near the same point $p$ in $e'(\br{F'})$ having 
maximal $z$ coordinate. Let $P$ be a horizontal plane passing
slightly below $p$
(so that in a neighborhood of $p$ it intersects $F$ in only one circle.)
We may pull the tube $e(A)$ until it is all above $P$. We then let it pass 
through itself until it is unknotted. This is a 0-extendible regular
homotopy with no quadruple points, at the end of which we have an embedding
intersecting $P$ as in Lemma \pr{s2} with an embedding of a torus
above the plane $P$, this embedding being already standard
and an embedding of a subsurface $F''$ of $F$
below the plane $P$, $F''$ being 
of smaller genus than that of $F$.
Again we are done by induction and Lemma \pr{s2}.

\emph{Case 4}: $D\su M^1(e)$ and $\pa D$ does not separate $e(F)$.
We may proceed as in Case 3 (this time via a 1-extendible regular homotopy)
to obtain a standard embedding of a torus connected with a 
tube to $e'(\br{F'})$
but this time the torus is contained in $M^0(e')$ and the tube connects
to $e'(\br{F'})$ from its compact side. But once we have such 
an embedding then the little standardly embedded
torus has non-separating compressing discs on 
both sides and so we are done by Case 3.

\end{pf}

\begin{lemma}\label{g5}

If $e:F\to\E$ is an embedding and $h:F\to F$ is a diffeomorphism such that
$e$ and $e\circ h$ are regularly homotopic, then
$\pp( A^0(e),A^1(e) ; A^0(e\circ h),A^1(e\circ h) ) = \p(h_*)$ and 
$\eep(e,e\circ h)=\ep(h)$.
(Recall that $h_*$ is the map induced by $h$ on $H_1(F,\C)$ and $\ep(h)\in\C$ 
is 0 or 1 according to whether $h$ is orientation
preserving or reversing.)
\end{lemma}

\begin{pf} 
$x\in \ker (e\circ h)^k_*$ iff $h_*(x) \in \ker e^k_*$
and so $A^k (e\circ h) = h_*^{-1}(A^k(e))$, $k=0,1$.
By definition then
$\pp( A^0(e),A^1(e) ; A^0(e\circ h),A^1(e\circ h) ) = \p(h_*^{-1})= \p(h_*)$.
(Note that if $e$ and $e\circ h$ are regularly homotopic then
indeed $h_*^{-1} \in O(\hh,g^e)$.)
$\eep(e,e\circ h)=\ep(h)$ is clear.
\end{pf}

We are now ready to prove Theorem \pr{t1}. For two regularly 
homotopic embeddings $e,e':F\to\E$ let  
$ \hp(e,e') = \pp(  A^0(e),A^1(e)  ;  A^0(e'),A^1(e')  )  
+  (n+1)\eep(e,e')$. We need to show $Q(e,e')=\hp(e,e')$.
If $e'':F\to \E$ is also in the same regular homotopy class
then $Q(e,e'')=Q(e,e')+Q(e',e'')$
and $\hp(e,e'')=\hp(e,e')+\hp(e',e'')$. And so if $e'$ is equivalent to
$e''$ and $Q(e,e'')=\hp(e,e'')$ then also $Q(e,e')=\hp(e,e')$.
And so we may replace $e$ with an equivalent standard embedding $f$
(Proposition \pr{p3})
and similarly replace
$e'$ with an equivalent standard embedding $f'$. Now $f$ and $f'$ 
have isotopic images and so after isotopy we may assume $f(F)=f'(F)$
and so $f'=f\circ h$ for some diffeomorphism $h:F\to F$.
By Lemma \pr{g5} and Theorem \pr{tn2} the proof of 
Theorem \pr{t1} is complete.

We conclude with a remark on systems of surfaces. 
If $S=F_1 \cup \cdots \cup F_r$ is a system of closed orientable surfaces, 
and $e:S\to \E$ is an
embedding, then we can rigidly move $e(F_i)$ one by one, until 
they are all contained in large disjoint balls. When it is the turn of
$F_i$ to be rigidly moved, then the union of all other components is embedded
and so only double lines occur. If $e':S\to \E$ is another embedding
then we can similarly move $e'(F_i)$ into the corresponding balls. 
It follows that $Q(e,e')=\sum_{i=1}^r Q(e|_{F_i},e'|_{F_i})$
and so we obtain a formula for systems of surfaces, namely:
$Q(e,e')=\sum_{i=1}^r \hp(e|_{F_i},e'|_{F_i})$.


\begin{thebibliography}{StoPC}

\bibitem[C]{C}
C.C. Chevalley: ``The Algebraic Theory of Spinors.'' 
Columbia University Press 1954. Also reprinted in:
C. Chevalley: ``The Algebraic Theory of Spinors and Clifford Algebras, 
Collected Works, Vol. 2'' 
Springer-Verlag 1997.

\bibitem[MB]{MB}
N. Max, T. Banchoff: ``Every Sphere Eversion Has a Quadruple Point.''
\emph{Contributions to Analysis and Geometry}, 
John Hopkins University Press,
(1981), 191-209. 

\bibitem[N1]{N1}
T. Nowik:
``Quadruple Points of Regular Homotopies of Surfaces in 3-Manifolds.''
\emph{Topology} - to appear.
(may be viewed at: 
http://www.math.columbia.edu/$\tilde{\ }$tahl/publications.html)

\bibitem[N2]{N2}
T. Nowik:
``Subgroups of the Mapping Class Group 
and Quadruple Points of Regular Homotopies.'' Preprint.
(may be viewed at: 
http://www.math.columbia.edu/$\tilde{\ }$tahl/publications.html)

\bibitem[P]{P}
U. Pinkall: ``Regular homotopy classes of immersed surfaces.''
\emph{Topology} 24 (1985) No.4, 421--434.


\end{thebibliography}
\end{document}